\documentclass{article}
\usepackage{amsmath}

\author{Gert Almkvist}
\title{Integrity of ghosts}
\begin{document}

\maketitle

\textbf{Introduction.}

Let $f$ be an $r\times r-$matrix with integer entries. Let further 
\[
\det (1+tf)=1+a_1t+a_2t^2+...+a_rt^r 
\]
be the characteristic polynomial of $f.$ Then $a_1,a_2,...a_r$ are of course
also integers. Conversely given integers $a_1,a_2,...,a_r$ we can find an $%
r\times r-matrix$%
\[
\left( 
\begin{array}{rrrrr}
0 & 0 & ... & 0 & \pm a_r \\ 
1 & 0 & ... & 0 & \mp a_{r-1} \\ 
... & ... & ... & ... & ... \\ 
0 & 0 & ... & 0 & -a_2 \\ 
0 & 0 & ... & 1 & a_1
\end{array}
\right) 
\]
with integer entries such that $\det (1+tf)=1+a_1t+...+a_rt^r$.

The \textsl{trace sequence }$b_1,b_2,...,b_r$ where $b_i=Tr(f^i)$ also
consists of integers. But these cannot be chosen arbitrarily. E.g. since $%
2a_2=b_1^2-b_2$ it follows that $b_1\equiv b_2$ mod $2$ . The aim of this
note is to find all congruences necessary and sufficient for $%
b_1,b_2,...,b_r $ to be a trace sequence. The condition is surprisingly
simple: 
\[
b_n\equiv b_{n/p}\text{ mod }p^k\text{ for }n\leq r 
\]
where $p$ is a prime such that 
\[
n=p^ks\text{ where }(p,s)=1 
\]
As a byproduct we find for for $g\in G$ that 
\[
\chi (g^{p^k})\equiv \chi (g^{p^{k-1}})\text{ mod }p^k 
\]
if $\chi $ is an integer valued (complex) character of a finite group (e-g.
the symmetric group or the Monster). This gives a partial answer to a
question by R.Brauer about conditions on a character table ( Problem 6,
p.139 in [3] ). The proof uses the fact that the trace sequence consists of
the ''ghost components'' of a certain Witt vector with integer coefficients.

The congruences up to $r=10$ were suggested by a computer program made by
Henrik Eriksson (in 1983) for which I am most thankful to him. 
\[
\]

\textbf{The proofs.}

We start with

\textbf{Proposition 1.}Given 
\[
\prod_{i=0}^\infty (1-x_it^i)=1-a_1t+a_2t^2-a_3t^3+a_4t^4-.... 
\]
then $x_n$ is a polynom with integer coefficients in $a_1,a_2,...,a_n.$
Similarly $a_n$ is a polynomial over \textbf{Z} in $x_1,x_2,...,x_n$%
\[
\]
$.$

\textbf{Proof: }$x_1=a_1.$The rest follows easily by induction. 
\[
\]

\textbf{Corollary. }All $a_i$ are integers if and only if all $x_i$ are
integers. 
\[
\]

The sequence 
\[
(x_1,x_2,...,x_n,....) 
\]
is called a \textsl{Witt vector. }For more about Witt vectors see [2]. 
\[
\]

\textbf{Definition 3. }The \textsl{ghost components }$b_1,b_2,...$ of the
Witt vector $(x_1,x_2,...)$ are given by 
\[
-t\frac d{dt}\log (\prod_{i=1}^\infty (1-x_it^i)=\sum_{n=1}^\infty b_nt^n 
\]
\[
\]

\textbf{Proposition 4. } 
\[
b_n=\sum_{d\mid n}dx_d^{n/d} 
\]
In particular if $p$ is a prime then 
\[
b_{p^k}=x_1^{p^k}+px_p^{p^{k-1}}+p^2x_{p^2}^{p^{k-2}}+...+p^kx_{p^k} 
\]
\[
\]

\textbf{Proof: }Expanding the left hand side of Definition 3 we get 
\[
\sum_{i=1}^\infty \sum_{j=0}^\infty ix_i^jt^{ij}=\sum_{n=1}^\infty \left\{
\sum_{i\mid n}ix_i^{n/i}\right\} t^n 
\]
\[
\]

\textbf{Theorem 5. Let }$f$ be a square matrix with integer entries. Then
for prime $p$ we have 
\[
Tr(f^{p^k})\equiv Tr(f^{p^{k-1}})\text{ mod}(p^k) 
\]
\[
\]

\textbf{Proof:} Let $\det (1-tf)=1-a_1t+a_2t^2-...\pm a_rt^r$. Then $%
a_{1,}a_2,...,a_r$ are integers and so are $x_1,x_2,...$ in the
corresponding Witt vector. Using the ''Exponential Trace Formula'' (see [1],
Theorem 1.10) 
\[
-t\frac d{dt}\log (\det (1-tf)=\sum_{i=1}^\infty Tr(f^i)t^i 
\]
we find that $b_n=Tr(f^n)$ for $n=1,2,..$are the ghosts of $x_1,x_2,...$. By
Proposition 4 we get 
\[
b_{p^k}=x_1^{p^k}+px_p^{p^{k-1}}+...+p^kx_{p^k} 
\]
\[
b_{p^{k-1}}=x_1^{p^{k-1}}+px_p^{p^{k-2}}+...+p^{k-1}x_{p^{k-1}} 
\]
It follows 
\[
b_{p^k}-b_{p^{k-1}}=\left\{ x_1^{p^k}-x_1^{p^{k-1}}\right\} +p\left\{
x_p^{p^{k-1}}-x_p^{p^{k-2}}\right\} +...+p^{k-1}\left\{
x_{p^{k-1}}^p-x_{p^{k-1}}\right\} +p^kx_{p^k} 
\]
We finish by showing the following Lemma. 
\[
\]

\textbf{Lemma 6.} For all integer $a$ we have 
\[
p^k\mid \left\{ a^{p^k}-a^{p^{k-1}}\right\} 
\]

\textbf{Proof.}We use induction over $k.$ For $k=1$ the result is Fermat's
little theorem. Assume that 
\[
a^{p^{k-1}}=a^{p^{k-2}}+cp^{k-1} 
\]
Now take the p$^{th}$ power of both sides.

\[
\]

\textbf{Corollary 7. } 
\[
Tr(f^{p^k})\equiv Tr(f^{p^{k-r}})\text{ mod }p^{k-r+1} 
\]
\[
\]

\textbf{Remark 8. }The Theorem is also valid when $f$ does not have integer
entries. We only need that $\det (1+tf)$ has integer coefficients. 
\[
\]

\textbf{Theorem 9. }Let $\chi $ be a complex character of a finite group $G$%
, that takes integer values. Then 
\[
\chi (g^{p^k})\equiv \chi (g^{p^{k-1}})\text{ mod }p^k 
\]
for all $g$ in $G.$ In particular it is true for all characters of the
symmetric group. 
\[
\]

\textbf{Proof. }Since $b_n=\chi (g^n)=Tr(g^n)\in \mathbf{Z}$ it follows by
the exponential trace formula (or by Newton's formulas for the symmetric
functions) that $a_i\in \mathbf{Q}$ . But $a_i$ is the the i$^{th}$
elementary symmetric function of the eigenvalues of $g$. Hence $a_i$ is a
rational algebraic integer so $a_i\in \mathbf{Z}$ .Then the rest follows
from Theorem 5. 
\[
\]

\textbf{Remark 10. }From Lemma 6 it follows that 
\[
\det (f^{p^k})\equiv \det (f^{p^{k-1}})\text{ mod }p^k 
\]
\[
\]
But we can do better than that. 
\[
\]

\textbf{Proposition 11.} Let $f$ be a square matrix with integer entries.
Then 
\[
\det (1+tf^{p^k})\equiv \det (1+tf^{p^k})\text{ mod p}^k 
\]
\[
\]

\textbf{Proof: }Let 
\[
\det (1+tf^{p^k})=1+c_1t+...+c_it^i+...+c_rt^r 
\]
Then 
\[
c_i=Tr(\Lambda ^i(f^{p^k}))=Tr((\Lambda ^if)^{p^k}) 
\]
since $\Lambda ^i$ is a functor. Hence by applying Theorem 5 with $f$
replaced by $\Lambda ^if$ we have 
\[
Tr((\Lambda ^if)^{p^k})\equiv Tr((\Lambda ^if)^{p^{k-1}})\text{ mod }p^k 
\]
\[
\]

\textbf{Theorem 13.} Let $b_1,b_2,...,b_r$ be integers. Then there exists an 
$r\times r-$matrix $f$ with integer entries such that $b_i=Tr(f^i)$ if and
only if the following condition is fulfilled for all $n\leq r:$

If $\quad n=mp^k$ where $(m,p)=1$ then 
\[
b_n\equiv b_{n/p}\text{ mod }p^k 
\]
\[
\]

\textbf{Proof:} The necessity of the condition follows from Theorem 5 with $%
f $ replaced by $f^m.$

To see the sufficiency we have to show that $x_n$ are integers for $n\leq r$
. Then it follows from Proposition 1 that $a_1,a_2,...,a_r$ are integers and
then we find $f$ as in the introduction.

We use induction on $n.$ If $n=1$ then $x_1=b_1$. Assume that $x_i\in 
\mathbf{Z}$ for $i<n.$ By Proposition 4 we have with $n=mp^k$ where $(m,p)=1$%
\[
b_n=b_{mp^k}=\sum_{d\mid m}\sum_{j=0}^kdp^j\left\{ x_{dp^j}^{m/d}\right\}
^{p^{k-j}} 
\]
We get 
\[
b_n-b_{n/p}= 
\]
\[
nx_n+\sum_{d\mid m,d<m}dp^kx_{dp^k}^{m/d}+\sum_{d\mid
m}\sum_{j=0}^{k-1}dp^j\left\{
(x_{dp^j}^{m/d})^{p^{k-j}}-(x_{dp^j}^{m/d})^{p^{k-j-1}}\right\} 
\]
By assumption $p^k\mid (b_n-b_{n/p})$. Using the induction hypothesis $%
x_{dp^j}$ is an integer for all $d\mid m$ and all $0\leq j<k.$ Lemma 6
implies that 
\[
p^{k-j}\mid \left\{
(x_{dp^j}^{m/d})^{p^{k-j}}-(x_{dp^j}^{m/d})^{p^{k-j-1}}\right\} 
\]
and hence every term except $nx_n$ on the right hand side is divisible by $%
p^k.$ It follows that 
\[
\frac{nx_n}{p^k} 
\]
is an integer. Do this with every prime $p$ dividing $n.$ It follows from
the Chinese Remainder Theorem that $x_n$ is an integer. 
\[
\]

\textbf{Ghosts in the Moonshine.}

It has been pointed out to me by John Walter that the results above could be
used on the monster group $M$ of Fischer and Griess. Following the notation
of Ogg ([8], p.531), for every $x$ in $M$ there is a modular function 
\[
J(x)=\frac 1q+\sum_{n=1}^\infty \omega _n(x)q^n 
\]
Here 
\[
q=\exp (2\pi ix) 
\]
and $\omega _n(x)$ are integers such that $x\rightarrow \omega _n(x)$ are
certain (reducible) characters of $M.$ ($J(x)$ depends only on the conjugate
class of $x$ in $M$. In particular $J(1)=j-744$ where 
\[
j=\frac 1q+744+1996884q+... 
\]
is the modular invariant (see [3], [7], [8], [11], [12] ). 
\[
\]

\textbf{Theorem. }Let $p$ be a prime. Then 
\[
J(x^{p^k})\equiv J(x^{p^{k-1}})\text{ mod }p^k 
\]
\[
\]

\textbf{Proof.} This follows immediately from Theorem 9 since the characters
take only integer values. 
\[
\]

\textbf{Remark.} The result is only interesting in case $p$ is one of the 15
primes ( $2\leq p\leq 31,$ $p=41,47,59,71$ ) dividing the order of $M$.
Since the maximal order of an element in $M$ is $119$, only a few small $k$
are of interest. 
\[
\]

\textbf{History of the problem}

The author ran into the problem of characterizing trace sequences in
connection with K-theory (see [1], p.296 ). After more than 10 years I found
the proof of Theorem 12, but first after Henrik Eriksson on a computer had
suggested the congruences up to $r=10.$

This paper was written in 1983. John McKay informed me that Frame [4] proved
Theorem 5 (for matrices of finite order ) already in 1949, using roots of
unity. After my talk at University of Illinois, Urbana in the Fall of 1983,
Stanley gave three more references: Carlitz 1958 [2], Schur 1937 [9] and
J\"{a}nischen 1921 [6]. My proof is almost identical with Schur's proof.
Thus he uses Witt vectors (probably before Witt). But he does not mention
matrices or any application to character tables. This is remarkable since he
generalizes some results to the integers of a number field. It is my hope
that putting this paper on the net shall make Schur's result more wellknown. 
\[
\]

\textbf{References.}

1. G. Almkvist, Endomorphisms of finitely generated projective modules over a
commutative ring, Arkiv f. Matematik, 11 (1973), 263-311.

2. L. Carlitz, Note on a paper of Dieudonn\'{e}, proc of AMS. 9 (1958), 32-33

3. J. Conway, S.Norton, Monstrous Moonshine, Bull. London Math. Soc. 11
(1979), 308-339.

4. J. S. Frame, Congruence relations between the traces of matrix powers,
Can.J.Math. 11 (1949), 303-304.

5. M. J. Greenberg, Lectures on forms in many variables, Benjamin, New York,
Amsterdam, 1959.

6. W. J\"{a}nischen, Ueber die Verallgemeinerung einer Gausschen Formel aus
der Theorie der h\"{o}heren Kongruenzen, Sitzungsberichte der Berliner Math.
Gesellschaft, 20 (1921), 23-29.

7. G. Mason, Modular forms and the theory of Thompson series,

8. A. Ogg, Modular functions, in Proc. of Symposia in Math. , Vol. 1, Wiley
1963

9. I. Schur, Aritmetische Eigenschaften der Potenzsummen einer algebraischen
Gleichung, Comp. math. 4 (1937), 432-444.

10. R. P. Stanley, Enumerative Combinatorics, Vol. 2, p.105, Problem 5.2,
Cambrodge Univ. Press 1999.

11. J. Thompson, Finite groups and modular functions, Bull. london math. Soc.
11 (1979), 347-351.

12. J. Thompson, Some numerology between the Fischer-Griess monster and the
elliptic modular function, Bull. London Math. Soc. 11 (1979), 352-353. 
\[
\]

Math Dept

Univ of Lund

Box 118

22100 Lund, Sweden

gert@maths.lth.se

\end{document}